\documentclass[11pt,a4paper,leqno]{amsart}
\usepackage {amsthm, amsmath, amssymb}
\usepackage[applemac]{inputenc}

\newtheorem{theorem}{Theorem}[section]

\theoremstyle{remark}
\newtheorem{remark}[theorem]{Remark}

\newcommand {\bfA}{\mathbf {A}}
\newcommand {\bfC} {\mathbf {C}}
\newcommand {\bfF} {\mathbf {F}}
\newcommand {\bfG} {\mathbf {G}}
\newcommand {\bfQ} {\mathbf {Q}}
\newcommand {\bfR} {\mathbf {R}}
\newcommand {\bfZ} {\mathbf {Z}}
\newcommand{\bfGL}{\mathbf{GL}}

\newcommand {\frakm} {\mathfrak{m}}

\newcommand {\numero} {n$^\circ$}
\newcommand {\Fbar} {\overline{\bfF}}

\DeclareMathOperator{\chr}{char}

\renewcommand* {\mod}{\mathop{\mathrm {mod}}\nolimits}

\title{How to use finite fields for problems concerning infinite fields}
\author{Jean-Pierre Serre}

\address{
Jean-Pierre Serre
\newline \indent
Coll\`ege de France
\newline \indent
3, rue d'Ulm, F-75005 Paris
}

\thanks{I want to thank A.Zykin who wrote a preliminary version of this lecture}
\date {\today}

\begin{document}
\maketitle

  As the title indicates, the purpose of the present lecture is  to show how to use finite fields for solving problems on infinite fields. This can be done on two different levels:  the elementary one uses only the fact that most algebraic geometry statements involve only finitely many data, hence come from geometry over a finitely generated ring, and  the residue fields of such a ring are finite; the examples we give in \S\S1-4 are of that type. A different level consists in using Chebotarev's density theorem and its variants, in order to obtain results over non-algebraically closed fields; we give such examples in \S\S5-6. The last two sections were only briefly mentioned in the actual lecture; they explain how cohomology (especially the \'etale one) can be used instead of  finite fields; the proofs are more sophisticated\footnote{Indeed, I would not have been able to give them without the help of Luc Illusie and of his two reports \cite{Ill} and \cite{Illusie}. }, but the results have a wider range.
  
  \section{Automorphisms of the affine $n$-space}

Let us start with the following  simple example:

\begin{theorem}
\label{involution}
Let $\sigma$ be an automorphism of the complex affine $n$-space $\bfC^n$, viewed as an algebraic variety. Assume that $\sigma^2=1.$ Then $\sigma$ has a fixed point.
\end{theorem}

Surprisingly enough this theorem can be proved by ``replacing $\bfC$ by a finite field". 

More generally:

\begin{theorem}
\label{fixedpoint}
Let $G$ be a finite $p$-group acting algebraically on the affine space $\bfA^n$ over an algebraically closed field $k$ with $\chr k \neq p$. Then the action of $G$ has a fixed point.
\end{theorem}
\noindent {\it Proof of Theorem 1.2}

 a) The case  $k=\Fbar_\ell$, where $\ell$ is a prime number $\neq p$

We may assume that the action of $G$ is defined over some finite extension $\bfF_{\ell^m}$ of  $\bfF_\ell.$ Then the group $G$ acts on the product $\bfF_{\ell^m} \times \dots \times \bfF_{\ell^m}.$ However, $G$ is a $p$-group and the number of elements of $\bfF_{\ell^m} \times \dots \times \bfF_{\ell^m}$ is not divisible by $p$. Hence there is an orbit consisting of one element, i.e. there is a fixed point for the action of $G.$

 b) Reduction to the case $k=\Fbar_\ell$

Since $G$ is finite, we can find a ring $\Lambda\subset \bfC$ finitely generated over $\bfZ,$ over which the action of $G$ can be defined. This means that the action of $G$ is given by  $$g(x_1,\dots, x_n)= (P_{g,1}(x_1, \dots, x_n), \dots, P_{g,n}(x_1, \dots, x_n)),$$ where the coefficients of the polynomials $P_{g,i}(x_1, \dots, x_n)$ belong to $\Lambda.$  Assume that there is no fixed point. The system of equations $$x_i - P_{g,i}(x_1, \dots, x_n)=0$$ has no solution in $\bfC$. Thus, by Hilbert's Nullstellensatz, there exist polynomials $Q_{g,i}(x_1, \dots, x_n)$ such that 
\begin{equation}
\label{Null}
\sum_{g, i}(x_i - P_{g,i}(x_1, \dots, x_n))Q_{g,i}(x_1, \dots, x_n)=1.
\end{equation}
By enlarging $\Lambda$ if necessary, we may assume that it contains $1/p$ and  the coefficients of  the $Q_{g,i}$'s. Let $\frakm$ be a maximal ideal of $\Lambda.$ Then the field $\Lambda/\frakm$ is finite (see e.g. \cite{Bourbaki}, p.$68$, cor.1), we have $\chr \Lambda/\frakm \neq p$ (since $p$ is invertible in $\Lambda$) and by (\ref{Null}) the conditions of the theorem hold for the algebraic closure of $\Lambda/\frakm.$ So we can apply part a) of the proof to get a contradiction.

\smallskip

\noindent
{\it Remark.} The technique of replacing a scheme $X$ of finite type over $k$ by a scheme over $\Lambda$ is sometimes called ``spreading out $X$"; its properties are described in  \cite{EGA4}, \S10.4.11 and \S17.9.7. 

\smallskip

\noindent
{\it Question.} Assume the hypotheses of Theorem 1.2. Let  $k_o$ be a subfield of $k$ such that the action of $G$ is defined over $k_o$. Does there exist a fixed point of  $G$  which is rational over $k_o$ ? Even the case $k = \bfC$, $k_o = \bfQ$, $ |G| = 2$, $n = 3$  does not seem to be known.

\smallskip

\noindent
{\small{\it Exercises} 

$1$. Let  $L$  be an infinite set of prime numbers. For every $p \in L$ , let  $k(p)$ be a denumerable field of characteristic $p$. Let  $A = \prod k(p)$ be the product of the $k(p)$'s.
Show that there exists a quotient of $A$ which is isomorphic to a subfield of $\bfC$.({\it Hint}. Use an ultrafilter on $L$.)

$2$. Let  $P_i(X_1,...,X_n)$ be a family of polynomials with coefficients in $\bfZ$. Show that the following properties are equivalent:

 a) The  $P_i$'s have a common zero in  $\bfC$.
 
 b) There exists an infinite set of primes  $p$  such that the  $P_i$'s have a common zero in $\bfF_p$.
 
 c) For every prime $p$, except a finite number, there exists a field of characteristic $p$ in which the $P_i$'s have a common zero.
 
 $3$. Assume the hypotheses of Theorem 1.2. Show that the number of fixed points of  $G$  is either infinite or $\equiv 1 \mod p$. ({\it Hint.} Suppose the set $S$ of fixed points is finite. Using the same argument as in the proof of Theorem 1.2, we may assume that the action of  $G$ is defined over a finite field $k_1$ with $q$ elements, with $(q,p) = 1$,  that the points of $S$ are rational over $k_1$, and that  $k_1$ contains the $p$-th roots of unity. We then get  $|S| \equiv q^n \mod p$, hence $|S| \equiv 1 \mod p$ since $q \equiv 1 \mod p$.) 
 
\noindent  [Smith's theory gives more: if  $S$ is finite, it has one element only, see \S7.4.]}
 
\section{Fixed points for finite group actions}
Consider a finite group $G$ of order $m$ acting on a $k$-variety $X$, where $k$ is  an algebraically closed field\footnote{The reader may interpret the word  ``$k$-variety" in the sense of FAC, i.e. as meaning a separated and reduced scheme of finite type over Spec($k$), cf. \cite{EGA4}, \S10.10. Since we are only interested in the $k$-points, the ``reduced" assumption has no importance.}. Let us assume  that there are a finite number of fixed points in $X(k)$ and that $G$ acts freely outside these points. 

\begin{theorem}
\label{fixed}
Suppose we have two actions of  $G$  on $X$ satisfying the above properties, the number of fixed points being $a$ and $a'$ respectively. Then $a \equiv a' \  (\mod m).$
\end{theorem}

\noindent {\it Sketch of  proof}

Assume first that  $X$, the two actions of $G$, and the fixed points, are defined over a finite field $\bfF_q$. We then have
$$
|X(\bfF_q)|\equiv a  \ (\mod m) \quad \text{and} \quad |X(\bfF_q)|\equiv a' \ (\mod m),
$$
hence $a  \equiv a' \ (\mod m)$.

  The general case is reduced to this case by an argument similar to (but less obvious than) the one given in \S1. [One replaces  $X$ by a separated scheme of finite type $X_0$ over a ring $\Lambda$ which is finitely generated over $\bfZ$, in such a way that the the two actions of $G$ extend to $X_0$; one also needs that the corresponding fixed points  $Y$  and  $Y'$  are finite and \'etale over Spec($\Lambda$), and that the two actions of $G$ on $X_0 - Y$ and on $X_0 - Y'$ are free. That all these conditions can be met is a consequence of  \cite{EGA4}, {\it loc.cit.}] One can then reduce modulo a maximal ideal of $\Lambda$.
  \smallskip

\noindent
{\it Remarks} 

 1) Theorem $1.1$ is a corollary of Theorem $2.1$. Indeed  the involution  $z\mapsto -z$ of $\bfC^n$ has only one fixed point, hence the number of fixed points of  any other involution  is either odd or infinite; it cannot be 0.

2) The results of $\S\S 1,2$ can also be proved by topological arguments, see $\S7.4$ below.
 
\section{Injectivity and surjectivity of maps between algebraic varieties}
The following theorem was proved independently by J.Ax and A.Grothen\-dieck in the $60$'s (see \cite{Ax}, \cite{GS}, p.184, and \cite{EGA4}, \S 10.4.11), and has been rediscovered several times since.

\begin{theorem}
Let $X$ be an algebraic variety over an algebraically closed field $k.$ If a morphism $f: X \to X$ is injective then it is bijective.
\end{theorem}

\noindent {\it Sketch of  proof}

Assume first that  $k$ is an algebraic closure of a finite field  $k_1$, and that  $f$  is defined over $k_1$. Then  $X(k_1)$ is finite; since  $$f : X(k_1) \to X(k_1)$$ is injective, it is bijective. The same argument, applied to the finite extensions of  $k_1$, shows that  $f$ : $X(k) \to X(k)$ is bijective. The case of an arbitrary algebraically closed field is reduced to the one above as  in \S\S1-2, by choosing a ring $\Lambda$ of finite type over $\bfZ$ over which $X$ and $f$ are defined, and reducing modulo a maximal ideal of $\Lambda$; for more details, see Grothendieck, \cite{EGA4}, {\it loc.cit.}

\smallskip

\noindent
{\it Remark.} When $k = \bfC$ a topological proof of this theorem was given by Borel in 1969 (see \cite{Borel}).

\section{ Nilpotent groups}

In 1955 M.Lazard (\cite{Lazard}) proved the following theorem:

\begin{theorem}
\label{nilp}
Let $G$ be an algebraic group over an algebraically closed field $k.$ If the underlying variety of $G$ is isomorphic to the affine space $\bfA^n$, for some  $n \geqslant 0$, then $G$ is nilpotent.
\end{theorem}
\noindent {\it Sketch of  proof}

  Lazard proves even more: he shows that  $G$  is {\it nilpotent of class} $\leqslant n$, i.e. that every iterated commutator of length $> n$ is equal to $1$. As in $\S1$, we may assume that  $k$ is an algebraic closure of a finite field $k_1$ and that  $G$  is defined over $k_1$. If $p =$ char $k$, the group  $G(k_1)$ is a finite $p$-group, hence is nilpotent. By applying this to the finite extensions of  $k_1$ one sees that  $G(k)$ is a 
 locally nilpotent group, i.e. is an increasing union of nilpotent groups. A further argument is needed to show that  $G(k)$ is indeed nilpotent of class $\leqslant n$, see \cite{Lazard}\footnote{Note the following misprints in \cite{Lazard}: in the last part of the proof of Lemma 1, ``$s  \geqslant t$" should be ``$s$ divisible by $t$"  and  ``$W = C_2(X,W)$" should be ``$W = C_2(X,V)$".}.

\smallskip

 Note that we may deduce Theorem \ref{nilp} from Theorem \ref{fixedpoint} together with the following standard result:

\begin{theorem}
\label{bor}
Let $G$ be a connected linear algebraic group. Then  either $G$ is nilpotent or it contains a one dimensional torus $\bfG_m$ as a subgroup. \end{theorem}

\noindent {\it Proof of Theorem \ref{nilp} from Theorem \ref{bor} and Theorem \ref{fixedpoint} } 

Assume that the underlying variety of $G$ is isomorphic to $\bfA^n$. In particular, $G$ is an affine variety; this is known to imply that the group $G$ can be embedded in some $\bfGL_N$, i.e. that $G$ is a linear group. Assume further that $G$ is not nilpotent; then $\bfG_m\subset G$ by Theorem \ref{bor}. Choose an element $\sigma$ of $\bfG_m$ of prime order $\ell\neq \chr k.$ The element $\sigma$ acts on $G\simeq \bfA^n$ by left translation and thus by Theorem \ref{fixedpoint} has a fixed point: contradiction!

\smallskip
\noindent
{\small{\it Exercise}. Give a topological proof of Theorem \ref{nilp} when $k = \bfC$ by using the fact that, for a non nilpotent connected group, either $H^1$ or $H^3$ is non zero. Extend this proof to arbitrary fields using $\ell$-adic cohomology.}

\section{Finite subgroups of $\bf GL$$_n(\bfQ)$}

The following theorem is a well-known result of Minkowski (see \cite{Mink}, \cite{Sch}, \cite{Serre}); it gives a multiplicative upper bound for the order of a finite subgroup of  $\bf GL$$_n(\bfQ)$.

\begin{theorem}
\label{Mink}
Let $\ell$ be a prime number. If $A$ is a subgroup of $\bf GL$$_n(\bfQ)$ of order $\ell^a$ then 
\begin{equation}
\label{bound}
a\leqslant M(n,\ell)=\left[\frac{n}{\ell - 1}\right ]+\left[\frac{n}{\ell(\ell - 1)}\right ]+\left[\frac{n}{\ell^2(\ell - 1)}\right ]+\dots.
\end{equation}
\end{theorem}

\noindent {\it Proof for $\ell \neq 2$} 

\
The idea is to reduce $\mod p$ for an appropriate choice of $p\neq \ell.$ First, we have $A\subset \bf GL$$_n(\bfZ[1/N])$ for a suitable $N \geqslant 1.$ If $p$ is sufficiently large ($p > 2$ is enough) and does not divide $N$, then the reduction $\mod p$ gives a subgroup $A'$  of $\bf GL$$_n(\bfF_p)$ such that $|A|=|A'|.$ Hence $|A|$ divides
 $ | \bf GL$$_n(\bfF_p)| $ , which is equal to $p^{n(n-1)/2}\prod_{i=1}^n(p^i-1).$ Let us choose $p$ such that its image in $(\bfZ/\ell^2\bfZ)^*$ is a generator of this group. This is always possible by Dirichlet's theorem on primes in arithmetic progressions since $(\bfZ/\ell^2\bfZ)^*$ is cyclic. Once $p$ is chosen in this way, $p^i-1$ is divisible by $\ell$ only if $i$ is divisible by $\ell - 1$ and in this case the $\ell$-adic valuation $v_{\ell}(p^i-1)$  of  $p^i-1$ is equal to $1+v_{\ell}(i).$  Hence 
   $v_\ell(|A|) \leqslant \sum (1 + v_\ell(i))$, where the sum is over the integers  $i$  with $1 \leqslant i \leqslant n$ which are divisible by $\ell - 1$; a simple computation shows that this sum is equal to $M(n,\ell)$ if $\ell \neq 2$, cf. \cite{Serre}, \S1.3.
\smallskip
\noindent
\begin{remark}
  In the case $\ell=2$ one has to replace the group $\bf GL$$_n$ by an orthogonal group $\bf O$$_n$ in order to get the desired bound.
  \end{remark}
\begin{remark} The result is optimal in the sense that, for every prime number $\ell$, there exist subgroups of $\bf GL$$_n(\bfQ)$ of order $\ell^{M(n,\ell)}$, see  \cite{Mink} and \cite{Serre}. These subgroups have the following Sylow-type property:
\end{remark}

\begin{theorem}
Let $A, P$ be two finite $\ell$-subgroups of  $\bf GL$$_n(\bfQ)$. Suppose that  $|P| =\ell^{M(n,\ell)}.$ Then $A$ is $\bf GL$$_n(\bfQ)$-conjugate to a subgroup of $P$.
\end{theorem}
[In particular, if  $|A| = |P|$,  then $A$ and $P$ are conjugate.]

\smallskip

\noindent {\it Sketch of  proof}

We only give the proof when  $\ell\neq 2$ (otherwise we have to consider orthogonal groups). Let us reduce $\mod \ p$ for a prime $p$ chosen as above. We get two $\ell$-subgroups $A$ and $P$ of $\bf GL$$_n(\bfF_p)$; by construction, $P$ is a Sylow subgroup of $\bf GL$$_n(\bfF_p)$. By Sylow's theorem, $A$ is conjugate to a subgroup of $P$, i.e. there exists an embedding $i : A \to B$ which is the restriction of an inner automorphism of $\bf GL$$_n(\bfF_p)$. The linear representations $A \to \bf GL$$_n(\bfQ)$ and $A\stackrel{i}{\rightarrow}P\to \bf GL$$_n(\bfQ)$ become isomorphic after reduction $\mod \ p.$ Since $p\neq \ell,$ a standard argument shows that they are isomorphic over $\bfQ_p$, hence also over $\bfQ$, and this completes the proof.

\noindent
\begin{remark}
In  \cite{Sch}, Schur gave another proof of Minkowski's Theorem using an  interesting lemma on characters of finite groups (\cite{Serre}, \S2.1, prop.1). He also extended  the theorem (by a different method) to arbitrary number fields (see \cite{Serre}, \S2.2).
\end{remark}

\section{Generalizations to other algebraic groups and fields}

A natural question is:  what happens in Theorem \ref{Mink} if  $\bfQ$ is replaced by an arbitrary field $k$ and  $\bf GL$$_n$  by an arbitrary reductive group $G$ ? This is answered in \cite{Serre}: roughly speaking, one can give a sharp bound for the order of a finite $\ell$-subgroup of $G(k)$ 
when one knows the root system of $G$  and the Galois group of the $\ell$-cyclotomic tower of  $k$ (one needs also to assume that  $G$ is ``of inner type", but this is automatic for the most interesting examples, such as $G_2, F_4, E_7$ or $E_8$). The proof follows Minkowski's method; the main difference is that Dirichlet's theorem on arithmetic progressions is replaced by a variant of the Chebotarev's density theorem which applies to every normal domain which is finitely generated over $\bfZ$ (see \cite{Serre65}, $\S2.7$ and  \cite{Serre}, $\S6.4$). As a sample, here is the case where  $k = \bfQ$ and  $G$  is of type  $E_8$:

\begin{theorem}
\label{E8}
Let $G$ be a group of type $E_8$ over $\bfQ$ and let $A$ be a finite subgroup of  $G(\bfQ)$. Then $ |A| $ divides the number\begin{equation*}
M(\bfQ, E_8)=2^{30}\cdot 3^{13}\cdot 5^5 \cdot 7^4 \cdot 11^2 \cdot 13^2 \cdot 19 \cdot 31.
\end{equation*}
\end{theorem}

\noindent {\it Sketch of  proof}

\
 We may assume that  $G$ is defined over $\bfZ[1/N]$ for some $N \geqslant 1$ and that  $A$ is an $\ell$-subgroup of $G(\bfZ[1/N])$  for some prime $\ell$ (the general bound is then obtained by multiplicativity). By reducing $\mod \ p$ for  $p$ large enough, we see that $|A|$ divides $|E_8(\bfF_p)|$. One knows (see e.g. \cite{Carter}, Theorem 9.4.10) that 
$$|E_8(\bfF_p)|=p^{120}(p^2-1)(p^8-1)(p^{12}-1)(p^{14}-1)(p^{18}-1)(p^{20}-1)(p^{24}-1)(p^{30}-1).$$
By choosing  $p$   as in  Minkowski's proof, one gets the desired bound. As an example, let us do the computation when   $\ell$ is equal to $3$:

 Choose  $p \equiv 2, 4, 5$ or $7 \ (\mod 9)$.  Then $p^2 - 1$ is divisible by 3 but not by 9. This implies that the 3-adic valuations of the eight factors  $p^2-1$, $p^8-1$, ..., $p^{30}-1$ are respectively: 1, 1, 2, 1, 3, 1, 2, 2. Their sum is 13, as claimed.

\noindent
\begin{remark}
\label{remE8}
The bound in Theorem \ref{E8} is optimal in the following sense: for every $\ell = 2,3,...,31$, there exists a group $G $ of type $E_8$ over $\bfQ$   such that  $G(\bfQ)$ contains a subgroup of order $2^{30}$, $3^{13}$, ..., 31, cf. \cite{Serre}, \S13.5. [Recall that there are three different groups of type $E_8$ over $\bfQ$, up to isomorphism; they are characterized by the structure of their $\bfR$-points. For most values of $\ell$, one can choose $G$ such that $G(\bfR)$ is compact.]
\end{remark}
\noindent
{\small{\it Exercise}. Let $K$ be a quadratic number field with discriminant $d$. Show that Theorem \ref{E8} is valid with $\bfQ$ replaced by $K$ and $M(\bfQ, E_8)$ replaced by $c(d)M(\bfQ, E_8)$,  where  $c(d)$ is defined by:

$c(8)= 2^8 $, $c(5)=5^5$, $c(13)=13^2$, $c(17)= 17^2$, $c(29) =29$, $c(37)= 37$, 

$c(41)= 41$, $c(61)= 61$ and $c(d)=1$ for the other values of $d$ (in particular those which are negative).}

\section{Proofs via Topology}
  As indicated at the end of  \S2, the results of \S\S1,2  can also be obtained - and sometimes improved - by using topological methods, based on cohomology (either standard, if the ground field is $\bfC$, or \'etale). There are several ways to do so; we  shall summarize a few  of them below. 
  
  The notation will be the following: $X$ is an algebraic variety over an algebraically closed field $k$ and  $G$ is a finite group which acts on $X$; the fixed point set of $G$  is denoted by $X^G$. We assume that $X$ is quasi-projective (in most applications, $X$ is affine), so that the quotient  variety $X/G$ is well defined. The cohomology groups $H^i(X)$  will be understood as the \'etale ones, with arbitrary support; in case we need cohomology with proper support, we shall write  $H_c^j(X)$, with a $c$ in index position. The letter $\ell$ will denote a prime number distinct from char $k$.
  \smallskip
  
 $7.1$.  {\it Using Cartan-Leray's spectral sequence}
 
   Suppose that  $G$ acts freely on  $X$. There is a Cartan-Leray spectral sequence (first defined in \cite {CL} in the context of standard sheaf cohomomology - for the case of étale cohomology, see \cite{GR}, \S8.5) $$H^i(G,H^j(X,C)) \Longrightarrow H^{i+j}(X/G,C), $$ where  $C$  is any finite abelian group. If  $X$  is the affine $n$-space $\bfA ^n$ and $|C|$ is prime to char $k$, then $H^j(X,C) = 0 $ for $j > 0$ and $H^0(X,C) = C$. In that case the spectral sequence degenerates and gives  $H^i(G,C) = H^i(X/G,C)$ for every $i$, i.e. $X/G$ has the same cohomology as the classifying space of $G$. Take now $C = \bfZ/\ell\bfZ$ and suppose that $\ell$ divides $|G|$. It is well known that $H^j(G,C)$ is non zero for infinitely many $j$'s, and that $H^j(X/G,C)$ is zero for $j > 2$.dim $X$: contradiction! 
   
   Conclusion: {\it the only finite groups which can act freely on $\bfA ^n$ are the $p$-groups, with} $p = $ char $k$. This gives another proof of Theorem 1.1.
   
   \smallskip
   {\small{\it Exercise}. Let $G$ be a finite $p$-group, with  $p = $ char $k$. Show that there exists a free action of $G$ on $\bfA ^n,$ provided that $n$ is large enough. ({\it Hint}: embed $G$ in a connected unipotent group.)}
   
   \smallskip
 $7.2.$ {\it Using Euler-Poincar\'e characteristics}
 
   Let  $\chi(X)$ be the Euler-Poincar\'e characteristic of $X$, relative to the $\ell$-adic cohomology. It is known (see \cite{Illusie}, \cite{Laumon}) that $\chi(X)$ does not depend on the choice of $\ell$, and that it coincides with the Euler-Poincar\'e characteristic of $X$ with proper support (Grothendieck-Laumon's theorem, cf.  \cite{Laumon} - in case $k = \bfC$, see the Appendix of \cite {KP}). In other words, we have 
   \smallskip
   
   $  \chi(X) = \sum (-1)^j \ $dim $H^j(X,\bfQ_\ell) = \sum (-1)^j \ $dim $H_c^j(X,\bfQ_\ell)$.
   
   \smallskip
   A useful property of $\chi$ is its {\it additivity} : if  $X$ is the disjoint union of locally closed subvarieties $X_{\lambda}$, then $\chi(X) = \sum \chi(X_{\lambda})$; this follows from the definition of $\chi(X)$ via cohomology with proper support.
   
   If $|G|$ is prime to char $k$, and $G$ acts freely on  $X$, one has  $\chi(X) = |G|.\chi(X/G)$, cf. \cite{Illusie}. (In particular $\chi(X) = 1$ implies $|G| = 1$; since $\chi(\bfA ^n) = 1$, we recover the statement given at the end of  $\S7.1$.) 
   
   Assume now that $G$ is an $\ell$-group. We have:
     $$ \chi(X^G) \equiv \chi(X) \ (\mod \ell).$$
  Indeed, using induction on $|G|$, one may assume that $G$ is cyclic of order $\ell$; in that case, $G$ acts freely on $Y = X - X^G$ and we have
  $$ \chi(X) = \chi(X^G) + \chi(Y) = \chi(X^G) + \ell.\chi(Y/G) \equiv \chi(X^G) \ (\mod \ell).$$
  In particular, {\it if $\chi(X)$ is not divisible by $\ell$, then $X^G$ is non empty.} This gives another proof of Theorem 1.2 (with the letter $p$ replaced by $\ell$).
   A similar argument applies to Theorem 2.1: with the notation of that theorem, one has  $a \equiv \chi(X)  \ (\mod m)$ and $a'  \equiv \chi(X)  \ (\mod m)$, hence $a \equiv a' \ (\mod m)$.
    
 \smallskip

 $7.3.$  {\it Using Lefschetz numbers} 
 
   If $s$ is an element of $G$, of order prime to char $k$, let  $t(s)$ be its Lefschetz number: 
   
   \smallskip
   
   $ t(s) = \sum (-1)^j.$Tr$(s, H^j(X,\bfQ_\ell)) = \sum (-1)^j.$Tr$(s, H_c^j(X,\bfQ_\ell))$, \ \
   cf. \cite{Illusie}.
   
   \smallskip

    It is known ({\it loc.cit.}) that $t(s) = 0$ if $s$ has no fixed point. If $X = \bfA ^n$, one has obviously  $t(s) = 1$. Hence {\it every automorphism of $\bfA ^n$, whose order is finite and prime to $\rm{char}($k$)$, has a fixed point}. In the special case where $k = \bfC$, this was proved by D.Petrie and J.D.Randall \cite{PR} by a similar argument, based on standard cohomology.
    
    \noindent {\it Remark}. For a general report on the possible actions of algebraic groups (not necessarily finite) on $\bfA ^n$, see H.Kraft \cite{Kraft}.

    \smallskip
 
  $7.4.$  {\it Using Smith theory}
  
    In the situation of Theorem 1.2, for $k = \bfC$,  Smith theory (cf. P.A.Smith \cite{Smith} and A.Borel \cite{BorelS}) gives more than the mere existence of a fixed point: it gives non trivial information on the cohomology of $X^G$. For instance, it shows that, if dim $X^G = 0$, then  $X^G$  is reduced to one point. Similar results hold in any characteristic. More precisely, suppose that $G$ is a finite $p$-group and let us write  $H^j(X)$ instead of $H^j(X,\bfZ/p\bfZ)$.
    
    \smallskip
      
 \noindent{\bf Theorem 7.5.} {\it Let $N$ be an integer such that $H^j(X) = 0$  for  all $j \geqslant  N$. Then} :
 
 \smallskip
   a)  $H^j(X^G) = 0$ {\it for all $j \geqslant  N$.
   
   \smallskip
    b) If $N = 1$ and  \rm{dim} $H^0(X) = 1$, {\it then} dim $H^0(X^G) = 1.$

 \smallskip
   The proof when $|G| = p$ will be given in \S8 below. The general case follows by induction on $|G|$: if $|G| > p$, one chooses a central subgroup $H$ of $G$ of order $p$ and one applies the induction hypothesis to the action of $G/H$ on $X^H$.
   
   \smallskip
 \noindent {\bf Corollary.} {\it If $X$ is $p$-acyclic, so is $X^G$}. 
 
 \smallskip
 This is obvious, since `` $p$-acyclic " means that dim $H^j$  is equal to 0 for $j > 0$ and to 1 for $j = 0$. Note that this implies that, if  $X^G$ is finite,  it has only one element.
 
 \noindent {\it Remark.} The prime number $p$ is allowed to be equal to char $k$. This
 case is of interest mainly when $X$ is a projective variety. For instance, if $X$ is a smooth projective surface in characteristic $p$ and if $X$ is rational, then  $X$ is $p$-acyclic [use the Artin-Schreier exact sequence] and the corollary above shows that the same is true for $X^G$; in particular, $X^G$ is not empty.

 \section{Smith Theory: Proof of Theorem 7.5 }
   Let $G$ be a cyclic group of prime order $p$. The group algebra  $\bfF_p[G]$ is isomorphic to the truncated polynomial ring $\bfF_p[t]/(t^p)$. This very simple fact is basic in Smith's proofs. It explains the rather artificial-looking definitions given below.
   \smallskip 
   
   \noindent 8.1. {\it $R$-abelian categories. Definitions}
   
     Let  $F$ be a field and let $n$ be an integer $  > 1.$ Let $R$ be the $F$-algebra generated by an element $t$ with the relation $t^n = 0$; it has for basis ${1,...,t^{n-1}}$. Let $C$ be an $R$-{\it abelian category}, i.e. an abelian category such that, for every pair of objects $A,B$ of $C$, Hom$_C(A,B)$ has an $R$-module structure, and the composition of maps is $R$-bilinear. In particular, $t$ defines an endomorphism $t_A$ of every object $A$ of  $C$, and we have $t_A^n$ = 0. If $i$ is a positive integer, the image of the endomorphism $t_A^i$  of  $A$  will be denoted by $t^iA$, and its kernel will be written $A_{t^i}$.  
     
     We have  $A \supset tA \supset t^2A \supset ... \supset t^nA = 0$, and  $A/A_{t^i} = t^iA$. 
     
     We shall say that  $A$ is {\it constant} if $tA = 0$, and that $A$ is {\it free} if the morphism $A/tA \rightarrow t^{n-1}A$  given by $t_A^{n-1}$ is an isomorphism, i.e. if $A_{t^{n-1}}= tA;$ this implies that all the quotients $t^iA/t^{i+1}A$ are isomorphic to $A/tA$ for  $ i = 1,...,n-1$.
      
\noindent {\it Example}. If $C$ is the category $Mod_R$ of all $R$-modules, the notion of freeness just defined coincides with the usual one. As for the ``constant" $R$-modules, they are the $F$-vector spaces with zero $t$-action.)
      
 \medskip

 \noindent 8.2. {\it The (I,A,B) setting}
 
   We now choose an exact sequence in $C$:
   
 \noindent $(8.2.1) \  \  0 \rightarrow I \rightarrow A \rightarrow B \rightarrow 0$
 
 \noindent such that $I$ is free and $B$ is constant.
 
 \smallskip
 
 \noindent {\bf Lemma 8.2.2.}  {\it For every  $i$  with $1\leqslant i \leqslant n-1$ we have} :
 
 \
\noindent (8.2.3)  {\it The natural map $A_{t^i}/t^{n-i}A \rightarrow B$ is an isomorphism.}
  
  \
\noindent  (8.2.4) {\it The natural map $A/t^iA \rightarrow B \oplus t^{n-i}A$ is an isomorphism.}

\noindent (In (8.2.4), the map  $A/t^iA \rightarrow  t^{n-i}A$ is induced by $t_A^{n-i}$.)
   
\noindent {\it Proof}. For every object  $Y$ of $C$, let us put  $h_i(Y) =  Y_{t^i}/t^{n-i}Y$. If 
$$ 0 \rightarrow Y' \rightarrow Y \rightarrow Y'' \rightarrow 0$$
is an exact sequence in $C$, we have an hexagonal exact sequence
$$\begin{array}{ccccccc}

&&h_i(Y')&\rightarrow&h_i(Y)&&\\
&\nearrow&&&&\searrow &\\
h_{n-i}(Y'')&&&&&&h_i(Y'')\\
&\nwarrow&&&&\swarrow &\\
&&h_{n-i}(Y)&\leftarrow&h_{n-i}(Y')&&
\end{array}
$$

Apply this to (8.2.1). Since $B$ is constant, we have $h_i(B) = B$, and since $I$ is free, we have $h_i(I) = 0$ for every $i$. We thus get the exact sequence
$$   0 = h_i(I) \ \longrightarrow  \; h_i(A) \ \longrightarrow \; h_i(B) = B \ \longrightarrow \; h_{n-i}(I) = 0,$$ which proves (8.2.3). As for (8.2.4), it follows from the fact that  $A/t^iA \rightarrow t^{n-i}A$ is surjective and that its kernel, by (8.2.3), is $h_i(A) = B$.

\smallskip
\noindent {\it Remark}. When $C = Mod_R$ , (8.2.1) implies that the $R$-module $A$ is a direct sum\footnote{possibly infinite : Smith theory does not require any finiteness assumption.} of indecomposable modules which are isomorphic to either $F$ or $R$ : in other words, the only ``Jordan blocks" which can occur are of rank either 1 or $n$.

\medskip
\noindent 8.3. {\it The main statements}

   We now consider another $R$-abelian category $C'$ and a cohomological functor $(H, \delta)$ on $C$ with values in $C'$ (cf. \cite{Tohoku}, \S2.1); we assume that this functor is compatible with the $R$-structures of $C$ and $C'$. For every $j \in \bfZ$, and every object $E$ of $C$, $H^j(E)$ is an object of $C'$; every exact sequence $0 \rightarrow E \rightarrow E' \rightarrow E''  \rightarrow 0$ in the category $C$ gives rise to an infinite exact sequence in $C'$:
   $$  ... \longrightarrow H^j(E) \longrightarrow H^j(E') \longrightarrow H^j(E'')  \stackrel{\delta_j} \longrightarrow H^{j+1}(E)  \longrightarrow ...$$
   We make the following assumptions:
   
   (8.3.1) $H^j(E) = 0$ for every $j < 0$ and every $E \in Ob(C)$. This implies that the functor $H^0$ is left exact.
   
   (8.3.2) For every $E \in Ob(C)$, one has $H^j(E) = 0$ for  $j$  large enough (which may depend on $E$).
   
   We are now going to apply the cohomological functor  $(H, \delta)$ to the $C$-objects $I,A,B$ of the exact sequence (8.2.1). The first result is:
   \smallskip
   
   \noindent {\bf Proposition 8.3.3.} {\it Let $N$ be a positive integer. Assume that $H^j(A) = 0 $ for all $j \geqslant N$. Then the same is true for $B$, i.e. $H^j(B) = 0$ for all $j \geqslant N$.}
   
   \smallskip
   In the special case $N = 1$ one can say more:
   \smallskip
   
   \noindent {\bf Proposition 8.3.4.} {\it Assume that $H^j(A) = 0$ for all $j > 0$, and that  $H^0(A)$ is constant, i.e. $t.H^0(A) = 0$. Then the same is true for $B$ and the natural map $H^0(A) \rightarrow H^0(B)$ is an isomorphism.}
 
 \medskip
 \noindent 8.4. {\it Proof of Proposition 8.3.3}
 \
 
  We prove first:

  \noindent {\bf Lemma 8.4.1}. {\it Let $m$ be a positive integer. Suppose that $H^j(A)$, $H^j(B)$, $H^j(t^iA)$ and $H^j(A_{t^i})$ are $0$ for  $ j = m+1$ and all $i = 1,...,n-1$. Suppose also that $H^m(A) = 0$. Then $H^m(B), H^m(t^iA)$ and $H^m(A_{t^i})$ are $0$.}
  
  \smallskip
 
 \noindent {\it Proof}. Since $H^{m+1}(t^iA) = 0$ the map $H^m(A) \rightarrow H^m(A/t^iA)$ is surjective, and this shows that $H^m(A/t^iA) = 0$. By (8.2.4) this implies that $H^m(B)$ and $H^m(t^{n-i}A)$ are 0, and since $n-i$ takes all values between 1 and $n-1$ we also have $H^m(t^iA) = 0$. By (8.2.3), $A_{t^i}$ is an extension of $B$ by $t^{n-i}A$; hence $H^m(A_{t^i})$ is also 0.
 \smallskip
 
  We can now prove Proposition 8.3.3. Indeed, by (8.3.2), we may choose an   $m > N$ such that the hypotheses of Lemma 8.4.1 are satisfied. One then uses descending induction on $m$. By Lemma 8.4.1, this is possible until we reach $m = N$; the proposition follows.  
  
 Note that we obtain at the same time the vanishing of  $H^j(t^iA)$ and $H^j(A_{t^i})$ for all $j \geqslant N$ and all $i = 1,...,n-1$.
 
 \medskip
 
 \noindent 8.5. {\it Proof of Proposition 8.3.4}
 
 \
   We apply Proposition 8.3.3 with $N = 1$. As mentioned above, we obtain the vanishing, not only of $H^1(B)$, but also of $H^1(tA)$ and $H^1(A_{t})$. Since $H^0$ is left exact, the sequence
   $$ 0 \longrightarrow H^0(A_t) \longrightarrow H^0(A) \stackrel{t} \longrightarrow H^0(A)$$
is exact, and since $t.H^0(A) = 0$ we see that  $H^0(A_t) \rightarrow H^0(A)$ is an isomorphism. Using the exact sequence 
$$ H^0(A_t) \longrightarrow H^0(A) \longrightarrow H^0(tA) \longrightarrow H^1(A_t) , $$
we deduce that $H^0(tA) = 0$. By (8.2.4) this implies that the natural map $H^0(A/t) \rightarrow H^0(B)$ is an isomorphism. But the map $H^0(A) \rightarrow H^0(A/tA)$ is also an isomorphism, since  $H^0(tA)$ and $H^1(tA)$ are both $0$.
Hence $H^0(A) \rightarrow H^0(B)$ is an isomorphism.

\medskip
\noindent 8.6. {\it Proof of Theorem 7.5}

  We can now prove Theorem 7.5 by applying what we have done to the case where $R =  \bfF_p[G] = \bfF_p[t]/(t^p)$ and $C$ is the category of the $ \bfF_p[G]$-\'etale sheaves over $X/G$, i.e. the sheaves which are killed by $p$, and are endowed with an action of $G$. 
  
  We take for $C'$ the categorie $Mod_R$ of all $R$-modules, and for cohomological functor the functor ``\'etale cohomology over $X/G$ "; hence, if $Y$ is a sheaf belonging to  $C$, $H^j(Y)$ is nothing else than $H^j(X/G,Y)$; condition (8.3.1) is obviously satisfied and condition (8.3.2) is a well known theorem of M.Artin, cf.  \cite{Artin}, \S 5.
  
  We take for sheaf  $A$  the direct image by the map $\pi : X \rightarrow X/G$ of the constant sheaf $\bfZ/p\bfZ$, i.e. $A = \pi_* (\bfZ/p\bfZ)$; we have $H^j(A) = H^j(X,\bfZ/p\bfZ) = H^j(X)$, cf. [11, \S5].  There is a natural action of $G$ on $A$, coming from the action of $G$ on the variety $X$.
  
   We take for sheaf $B$ the constant sheaf $\bfZ/p\bfZ$ on $X^G$ ``extended by zero on $ X/G - X^G$ "  (direct image by the inclusion\footnote{ \ The natural map $X^G \rightarrow \pi(X^G)$ is a homeomorphism for the étale topology (it is even an isomorphism
   of schemes if $p$ is distinct from char $k$); hence we may identify $X^G$ with its image in $X/G$.}$X^G \rightarrow X/G$); we have $H^j(B) = H^j(X^G)$. The sheaf $B$ is constant (in the sense of 8.1) : the group $G$ acts trivially on it.
   
   There is a natural surjection $A \rightarrow B$. Its kernel $I$ is free in the sense of $\S8.1$. This is checked "fiber by fiber"; if $x$ is a geometric point of $X/G$, and $x$ belongs to $X^G$, the fiber of $I$ at $x$ is 0; if $x$ belongs to $X - X^G$, the fiber
   of $I$ at $x$ is a free $R$-module of rank 1.
   
   We can now apply Proposition 8.3.3 and Proposition 8.3.4. They give the two parts of Theorem 7.5.
   
   \noindent {\it Remark}. The same proof applies in the usual context of sheaf theory, provided that the space $X/G$ has finite cohomogical dimension (so that (8.3.2) holds). It also applies in a combinatorial context, with cohomology replaced by homology; this was already mentioned in  \cite{SerreB}, proof of Prop.2.10.

\end{document}